\documentclass[reqno,12pt,a4paper]{amsart}

\usepackage{mathrsfs}
\usepackage{amssymb}
\usepackage{amsfonts}
\usepackage{latexsym}
\usepackage{amsthm}

\usepackage{graphicx}
\def\lb{\label}

\newcommand{\er}[1]{\textrm{(\ref{#1})}}

\begin{document}


\renewcommand{\theequation}{\arabic{section}.\arabic{equation}}
\theoremstyle{plain}
\newtheorem{theorem}{\bf Theorem}[section]
\newtheorem{lemma}[theorem]{\bf Lemma}
\newtheorem{corollary}[theorem]{\bf Corollary}
\newtheorem{proposition}[theorem]{\bf Proposition}
\newtheorem{definition}[theorem]{\bf Definition}
\newtheorem{remark}[theorem]{\it Remark}

\def\a{\alpha}  \def\cA{{\mathcal A}}     \def\bA{{\bf A}}  \def\mA{{\mathscr A}}
\def\b{\beta}   \def\cB{{\mathcal B}}     \def\bB{{\bf B}}  \def\mB{{\mathscr B}}
\def\g{\gamma}  \def\cC{{\mathcal C}}     \def\bC{{\bf C}}  \def\mC{{\mathscr C}}
\def\G{\Gamma}  \def\cD{{\mathcal D}}     \def\bD{{\bf D}}  \def\mD{{\mathscr D}}
\def\d{\delta}  \def\cE{{\mathcal E}}     \def\bE{{\bf E}}  \def\mE{{\mathscr E}}
\def\D{\Delta}  \def\cF{{\mathcal F}}     \def\bF{{\bf F}}  \def\mF{{\mathscr F}}
\def\c{\chi}    \def\cG{{\mathcal G}}     \def\bG{{\bf G}}  \def\mG{{\mathscr G}}
\def\z{\zeta}   \def\cH{{\mathcal H}}     \def\bH{{\bf H}}  \def\mH{{\mathscr H}}
\def\e{\eta}    \def\cI{{\mathcal I}}     \def\bI{{\bf I}}  \def\mI{{\mathscr I}}
\def\p{\psi}    \def\cJ{{\mathcal J}}     \def\bJ{{\bf J}}  \def\mJ{{\mathscr J}}
\def\vT{\Theta} \def\cK{{\mathcal K}}     \def\bK{{\bf K}}  \def\mK{{\mathscr K}}
\def\k{\kappa}  \def\cL{{\mathcal L}}     \def\bL{{\bf L}}  \def\mL{{\mathscr L}}
\def\l{\lambda} \def\cM{{\mathcal M}}     \def\bM{{\bf M}}  \def\mM{{\mathscr M}}
\def\L{\Lambda} \def\cN{{\mathcal N}}     \def\bN{{\bf N}}  \def\mN{{\mathscr N}}
\def\m{\mu}     \def\cO{{\mathcal O}}     \def\bO{{\bf O}}  \def\mO{{\mathscr O}}
\def\n{\nu}     \def\cP{{\mathcal P}}     \def\bP{{\bf P}}  \def\mP{{\mathscr P}}
\def\r{\rho}    \def\cQ{{\mathcal Q}}     \def\bQ{{\bf Q}}  \def\mQ{{\mathscr Q}}
\def\s{\sigma}  \def\cR{{\mathcal R}}     \def\bR{{\bf R}}  \def\mR{{\mathscr R}}
\def\S{\Sigma}  \def\cS{{\mathcal S}}     \def\bS{{\bf S}}  \def\mS{{\mathscr S}}
\def\t{\tau}    \def\cT{{\mathcal T}}     \def\bT{{\bf T}}  \def\mT{{\mathscr T}}
\def\f{\phi}    \def\cU{{\mathcal U}}     \def\bU{{\bf U}}  \def\mU{{\mathscr U}}
\def\F{\Phi}    \def\cV{{\mathcal V}}     \def\bV{{\bf V}}  \def\mV{{\mathscr V}}
\def\P{\Psi}    \def\cW{{\mathcal W}}     \def\bW{{\bf W}}  \def\mW{{\mathscr W}}
\def\o{\omega}  \def\cX{{\mathcal X}}     \def\bX{{\bf X}}  \def\mX{{\mathscr X}}
\def\x{\xi}     \def\cY{{\mathcal Y}}     \def\bY{{\bf Y}}  \def\mY{{\mathscr Y}}
\def\X{\Xi}     \def\cZ{{\mathcal Z}}     \def\bZ{{\bf Z}}  \def\mZ{{\mathscr Z}}
\def\O{\Omega}

\newcommand{\gA}{\mathfrak{A}}
\newcommand{\gB}{\mathfrak{B}}
\newcommand{\gC}{\mathfrak{C}}
\newcommand{\gD}{\mathfrak{D}}
\newcommand{\gE}{\mathfrak{E}}
\newcommand{\gF}{\mathfrak{F}}
\newcommand{\gG}{\mathfrak{G}}
\newcommand{\gH}{\mathfrak{H}}   \newcommand{\gh}{\mathfrak{h}}
\newcommand{\gI}{\mathfrak{I}}
\newcommand{\gJ}{\mathfrak{J}}
\newcommand{\gK}{\mathfrak{K}}
\newcommand{\gL}{\mathfrak{L}}
\newcommand{\gM}{\mathfrak{M}}
\newcommand{\gN}{\mathfrak{N}}
\newcommand{\gO}{\mathfrak{O}}
\newcommand{\gP}{\mathfrak{P}}
\newcommand{\gQ}{\mathfrak{Q}}
\newcommand{\gR}{\mathfrak{R}}
\newcommand{\gS}{\mathfrak{S}}
\newcommand{\gT}{\mathfrak{T}}
\newcommand{\gU}{\mathfrak{U}}
\newcommand{\gV}{\mathfrak{V}}
\newcommand{\gW}{\mathfrak{W}}
\newcommand{\gX}{\mathfrak{X}}
\newcommand{\gY}{\mathfrak{Y}}
\newcommand{\gZ}{\mathfrak{Z}}

\def\ve{\varepsilon}   \def\vt{\vartheta}    \def\vp{\varphi}    \def\vk{\varkappa}

\def\Z{{\mathbb Z}}    \def\R{{\mathbb R}}   \def\C{{\mathbb C}}    \def\K{{\mathbb K}}
\def\T{{\mathbb T}}    \def\N{{\mathbb N}}   \def\dD{{\mathbb D}}


\def\la{\leftarrow}              \def\ra{\rightarrow}            \def\Ra{\Rightarrow}
\def\ua{\uparrow}                \def\da{\downarrow}
\def\lra{\leftrightarrow}        \def\Lra{\Leftrightarrow}


\def\lt{\biggl}                  \def\rt{\biggr}
\def\ol{\overline}               \def\wt{\widetilde}
\def\no{\noindent}


\let\ge\geqslant                 \let\le\leqslant
\def\lan{\langle}                \def\ran{\rangle}
\def\/{\over}                    \def\iy{\infty}
\def\sm{\setminus}               \def\es{\emptyset}
\def\ss{\subset}                 \def\ts{\times}
\def\pa{\partial}                \def\os{\oplus}
\def\om{\ominus}                 \def\ev{\equiv}
\def\iint{\int\!\!\!\int}        \def\iintt{\mathop{\int\!\!\int\!\!\dots\!\!\int}\limits}
\def\el2{\ell^{\,2}}             \def\1{1\!\!1}
\def\sh{\sharp}
\def\wh{\widehat}
\def\bs{\backslash}

\def\sh{\mathop{\mathrm{sh}}\nolimits}
\def\ch{\mathop{\mathrm{ch}}\nolimits}
\def\where{\mathop{\mathrm{where}}\nolimits}
\def\all{\mathop{\mathrm{all}}\nolimits}
\def\as{\mathop{\mathrm{as}}\nolimits}
\def\Area{\mathop{\mathrm{Area}}\nolimits}
\def\arg{\mathop{\mathrm{arg}}\nolimits}
\def\const{\mathop{\mathrm{const}}\nolimits}
\def\det{\mathop{\mathrm{det}}\nolimits}
\def\diag{\mathop{\mathrm{diag}}\nolimits}
\def\diam{\mathop{\mathrm{diam}}\nolimits}
\def\dim{\mathop{\mathrm{dim}}\nolimits}
\def\dist{\mathop{\mathrm{dist}}\nolimits}
\def\Im{\mathop{\mathrm{Im}}\nolimits}
\def\Iso{\mathop{\mathrm{Iso}}\nolimits}
\def\Ker{\mathop{\mathrm{Ker}}\nolimits}
\def\Lip{\mathop{\mathrm{Lip}}\nolimits}
\def\rank{\mathop{\mathrm{rank}}\limits}
\def\Ran{\mathop{\mathrm{Ran}}\nolimits}
\def\Re{\mathop{\mathrm{Re}}\nolimits}
\def\Res{\mathop{\mathrm{Res}}\nolimits}
\def\res{\mathop{\mathrm{res}}\limits}
\def\sign{\mathop{\mathrm{sign}}\nolimits}
\def\span{\mathop{\mathrm{span}}\nolimits}
\def\supp{\mathop{\mathrm{supp}}\nolimits}
\def\Tr{\mathop{\mathrm{Tr}}\nolimits}
\def\BBox{\hspace{1mm}\vrule height6pt width5.5pt depth0pt \hspace{6pt}}


\newcommand\nh[2]{\widehat{#1}\vphantom{#1}^{(#2)}}
\def\dia{\diamond}

\def\Oplus{\bigoplus\nolimits}



\def\qqq{\qquad}
\def\qq{\quad}
\let\ge\geqslant
\let\le\leqslant
\let\geq\geqslant
\let\leq\leqslant
\newcommand{\ca}{\begin{cases}}
\newcommand{\ac}{\end{cases}}
\newcommand{\ma}{\begin{pmatrix}}
\newcommand{\am}{\end{pmatrix}}
\renewcommand{\[}{\begin{equation}}
\renewcommand{\]}{\end{equation}}
\def\eq{\begin{equation}}
\def\qe{\end{equation}}
\def\[{\begin{equation}}
\def\bu{\bullet}

\title[{ KdV Hamiltonian as  function  of actions }]
{KdV Hamiltonian as  function  of actions}

\date{\today}
\author[Evgeny L. Korotyaev]{Evgeny L. Korotyaev}
\address{
Saint-Petersburg, Russia.
 \ korotyaev@gmail.com}
\author[Sergei Kuksin]{Sergei Kuksin}
\address{CMLS, Ecole Polytechnique, 91128 Palaiseau,
France,   kuksin@math.polytechnique.fr}


\begin{abstract}
We prove that the non-linear part of the Hamiltonian of the KdV  equation on the circle,
written as a function of the actions, defines a continuous convex function on the $\ell^2$ space and derive for it
lower and upper bounds in terms of some functions of the  $\ell^2$-norm.
The proof is based on a new representation of the Hamiltonian in terms of
the quasimomentum and its analysis using the conformal mapping theory.
\end{abstract}

\maketitle


\section {Introduction and Main results }
\setcounter{equation}{0}

We consider the Korteweg de Vries (KdV)  equation under zero mean-value periodic boundary conditions:
\[
\begin{aligned}
\lb{KdV}
& q_t=-q_{xxx}+6qq_x, \qqq  \qqq \qq x\in \T=\R/\Z,
\\
& \int _0^1q(x,t)\,dx=0.
\end{aligned}
\]
For any $\a\in \R$  denote by $\mH_\a$ the Sobolev space  of real-valued 1-periodic functions
with zero mean-value. In particular, we have
$$
\mH_\a=\mH_\a(\T)=\rt\{q\in L^2(\T): \ q^{(\a)}\in L^2(\T),\qq \int
_0^1q(x)\,dx=0\rt\},\qqq \a\ge 0.
$$
We provide the spaces $\mH_\a$ with the trigonometric base $\{e_{\pm 1}, e_{\pm 2}, e_{\pm 3}, ....\}$, where
$$
e_j(x)=\sqrt 2\cos 2\pi j x, \qqq e_{-j}(x)=-\sqrt 2\sin 2\pi j x,\qqq j\ge 1.
$$
We also introduce real spaces $\ell^p_{\a}$ of sequences $f=(f_n)_1^{\iy }$, equipped with the norms
\[
\lb{1.0}
\|f\|_{p,\a}^p=\sum _{n\ge 1}(2\pi n)^{2\a}|f_n|^p,\qqq p\ge 1,\ \a\in \R,
\]
and positive octants
$$
\ell^p_{\a, +}=\rt\{f=(f_n)_1^{\iy } \in \ell^p_{\a}: f_n\ge 0, \ \forall \ n\ge 1\rt\}.
$$
In the case $\a=0$ we write $\ell^p=\ell_0^p, \ell_+^p=\ell_{0,+}^p$ and
$\|\cdot\|_{p}=\|\cdot\|_{p,0}$.

The operator ${\pa \/\pa x}$ defines  linear isomorphisms ${\pa \/\pa x}: \mH_\a\to \mH_{\a-1}$.
Denoting by $({\pa \/\pa x})^{-1}$ the inverse operator, we provide the spaces $\mH_\a, \a\ge0$, with a symplectic
structure by means of the 2-form $\o_2$:
$$
\o_2(q_1,q_2)=-\lan ({\pa /\pa x})^{-1}q_1,q_2\ran,
$$
where $\lan \cdot,\cdot\ran$ is the scalar product in $L^2(0,1)$.
Then in any space $\mH_\a, \a\ge 1$, the KdV equation \er{KdV}  may be written as a Hamiltonian system
with the   Hamiltonian $H_2$, given by
$$
H_2(q)={1\/2}\int_0^1(q'(x)^2+2q^3(x))\,dx.
$$
That is,  as the  system
\[
\lb{KdV1}
q_t={\pa \/\pa x}{\pa \/\pa q}H_2(q),
\]
e.g. see \cite{Ku, KP}
(note that $H_2$ is an analytic function on any space $\mH_\alpha,\  \alpha\ge1$).

 It is well known  after the celebrated work of Novikov, Lax, Its and Matveev that   the system \er{KdV1} is integrable.
It was shown by Kappeler and collaborators in a series of publications, starting with \cite{Ka},
that it admits global Birkhoff coordinates. Namely, for any $\a\in \R$ denote by
$\gh_\a$ the Hilbert space, formed by real sequences $b=(b_n, b_{-n})_1^\iy$, equipped with the norm
$$
\|b\|_{\gh_\a}^2=\sum_{n\ge 1}(2\pi n)^{2\a}(b_n^2+b_{-n}^2).
$$
We provide  the spaces
$\gh_\a$ with the usual symplectic form
$$
\O_2=\sum_{n\ge 1}db_n \wedge db_{-n},
$$
and define the actions $I=(I_n)_1^\iy$ and the angles $\f=(\f_n)_1^\iy$ by
\[
\lb{1.22}
I_n={1\/2}(b_n^2+b_{-n}^2),\qqq \qqq \f_n=\arctan {b_n\/b_{-n}}.
\]
This is another set of symplectic coordinates on $\gh_\a$, since $\O_2=dI \wedge d\f$, at least formally. Then

{\it 1) There exists an analytic  symplectomorphism $\P: \mH_0\to \gh_{1\/2}$ which defines
analytic  diffeomorphisms
$\P: \mH_\a\to \gh_{\a+{1\/2}}, \a\ge -1$, such that $d\P(0)=\F$, where
\[
\lb{1x} \F\rt(\sum _{j\ge 1}\rt(u_je_j(x)+u_{-j}e_{-j}(x)\rt)
\rt)=b,\qqq b_j=|2\pi j|^{1\/2}u_j,\  \ \forall \ j.
\]

2) The transformed Hamiltonian $H_2(\P^{-1}(b))$ (which is an anlytic function on the space $\gh_{{3\/2}}$)
depends solely on the actions $I$, i.e. $K(I(b))=H_2(\P^{-1}(b))$,
 where $K(I)$ is an analytic function on the octant $\ell_{{3\/2},+}^1$.
 A curve $q(\cdot,t)\in C^1(\R, \mH_0)$ is a solution of \er{KdV} if and only if
$b(t)=\P(q(\cdot,t))$ satisfies the following system of equations
\[
\lb{hameq}
{\pa b_n\/\pa t}=-b_{-n}{\pa K\/\pa I_n},\qqq \qqq {\pa b_{-n}\/\pa t}=b_{n}{\pa K\/\pa I_n},
\qqq n\ge 1,
\]
where $I=(I(b))$.

}
\medskip

For 1)-2) with $\a\ge 0$  see \cite{KP} and with $\a=-1$ see \cite{KT}. See
\cite{KuP} for the important  quasilinearity property of the transformation $\Psi$.

Note that $q\in \mH_{-1}$ iff $I\in \ell_{-{1\/2}}^1$. Thus  if $I\in \ell^p$ for some $p\in [1,\iy)$,
then $I\in \ell_{-{1\/2}}^1$ and the corresponding potential $q\in \mH_{-1}$.

\smallskip

By 2), in the action-angle variables $(I,\phi)$ the KdV equation takes the form
\begin{equation}
\label{Ac_An}
I_t=0,\qquad \phi_t =\frac{\partial}{\partial I}K(I).
\end{equation}
This reduction of KdV
 is due to McKean-Trubowitz \cite{MT1} and was found before the Birkhoff form \eqref{hameq}.
The action maps $\psi\mapsto I_j$, $j\ge1$, are given by explicit formulas due to Arnold and are defined in
a unique way. So the Hamiltonian $K(I)$ also is uniquelly defined, see \cite{FM}. But the symplectic
angles are defined only up to rotations $\phi\mapsto  \phi+({\partial}/{\partial I})g(I)$, where $g$ is any
smooth function. So the transformation $\Psi$ is not unique.
\medskip

The Birkhoff coordinates $b$  and the actions-angles $(I,\phi)$
make an effective tool to study properties of the KdV equation, see \cite{KT},
and of  its perturbations, see \cite{Ku1}. For both these goals it is important to understand properties of the
Hamiltonian $K(I)$ which defines  the dynamics \er{hameq} and \eqref{Ac_An}.

Denote by $P_j$  moments of the actions $I$, given by
\[
\lb{1.70}
P_j=\sum_{n\ge 1} (2\pi n)^{j}I_n, \qqq \qqq \qq j\in\Z.
\]
Note that
\[
\lb{pq}
P_1={1\/2}\|q\|^2, \qqq \text{if}\qq I=I(b), \ b=\P(q),
\]
-- this is the Parseval identity for the transformation  $\P$, see  \cite{MT, K5}.
 Due to \er{1x}, the linear part $dK(0)(I)$ of $K(I)=H_2(\P^{-1}(b))$ at the origin  equals
 $$
 {1\/2}\int_0^1\rt( {\pa \/\pa x}\rt( \F^{-1} b\rt)\rt)^2\,dx=P_3.
 $$
Therefore
 $$
 K(I)=P_3(I)+O(I^2).
 $$
The cubic part $\int_0^1q^3(x)\,dx$ of the Hamiltonian $H_2(q)$ is more regular than its quadratic
part ${1\/2}\int_0^1q'(x)^2\,dx$. Thus  it is natural to assume that the term $P_3$
is a singular part of $K(I)$ and to study  smoothness of the more regular  quadratic
part $V$, given by
\[
\lb{HVid}
H_2(q)=K(I)=P_3(I)-V(I).
\]
Here the minus-sign is convenient, since below we will see that $V\ge 0$.
For any $N\ge 1$ denote by $\wt\ell^N\ss\ell^2$ the N-dimensional subspace
$$
\wt\ell^N=\{h=(h_n)_1^\iy, \ h_n=0\;\forall  \ n>N\}, \qqq
$$
and set $\wt\ell^\iy=\cup \wt\ell^N$.
Clearly $V$ is analytic on each octant $\wt\ell_+^N$ (i.e., it analytically extends to a neighbourhood
 of $\wt\ell_+^N$ in $\wt\ell^N$). So $V$ is Gato-analytic on $\wt\ell_+^\iy$. That is, it is analytic on each interval
 $\{  (a+tc)\in \wt\ell_+^\iy| t\in \R  \}$, where $a,c\in \wt\ell_+^\iy$. It is known that
\[
\lb{01}
 {\pa^2 V(0)\/\pa I_i\pa I_j}=6\d_{i,j}\qqq \forall \qq i,j\ge 1,
\]
see \cite{BoK} and \cite{KP, Ku}. So $d^2V(0)(I)=6\|I\|_2^2$.
This suggests that the Hilbert space $\ell^2$ rather than the Banach space $\ell^1_{3\/2}$
(which is contained in $\ell^{2}$)  is a distinguished phase-space for the
Hamiltonian $K(I)$. This guess is justified by the following theorem which is the main result of our work.

\

\begin{theorem}
\lb{T1}
The function $V: \wt\ell^\iy_+ \to \R$ extends to a non-negative continuous function on the $\ell^{2}$-octant $\ell_+^{2}$, such that $V(I)=0$ for some $I\in\ell_+^{2}$ iff $I=0$. Moreover,
\[
\lb{VA}
0\le V(I)\le 8P_1P_{-1},\qqq \forall \ I\in \ell_{1,+}^1,
\]
and
\[
\lb{VA2} {\pi\/10} { \|I\|_2^2\/(1+P_{-1}^{1\/2})}\le V\le \rt(4^{11\/2}(1+P_{-1}^{1\/2})^{1\/2}P_{-1}^2+ 6\pi
e^{\sqrt{P_{-1}}}\|I\|_2\rt)\|I\|_2, \qqq \forall \ I\in \ell^2.
\]
\end{theorem}

Let $X$ be a Banach space which contains $\tilde\ell^\infty$ as a dense subsets.
We say that the function $V(I)$ {\it agrees with the norm} $\|I\|_X$ if $V$ extends to a continuous function
on $X_+$ (= the closure of  $\wt\ell_+^\iy$ in $X$) and
$$
F_1(\|I\|_X)\le V(I)\le F_2(\|I\|_X),\qqq \qq \forall \ I\in X_+,
$$
where $F_1, F_2$ are monotonous continuous functions from $\R_+$ into $\R_+$ such that $F_j(0)=0$ and
$F_j(t)\to \iy$ as $t\to \iy$, $ j=1,2$.
It is easy to see that there exists at most one Banach space $X$ as above (i.e., if $X'$ is
another space, then $X=X'$ and the two norms are equivalent).

Estimates \er{VA2} imply that  the function  $V(I)$ agrees with the norm  $\| I\|_2$.
So $\ell^{2}$ is the natural phase space  for the non-linear part $V$ of the  Hamiltonian $K(I)$.
Estimates from Section 3 easily imply that $\ell^2$-sequences $I$ correspond to
potentials $q\in \mH_{-1}$ and in general these potentials do not belong to $\mH_{1/2}$ (see Remark 2 in Section 3).

A proof of the theorem  is based on a new identity (see Theorem \ref{TVe}), representing $V(I)$ in terms of the quasimomentum of the Hill operator with a potential $q$. It uses properties of the conformal mapping, associated
with this quasimomentum,  developed in \cite{K1} - \cite{K5}.

{\bf Remarks.} 1) \er{VA} improves the known estimate
$|H_2(q)|\le 4^5P_3(1+P_3^{4\/3})$ from \cite{K3}.

2) We claim that the function $V$ is real analytic on $\ell_+^2$. This will be proven elsewhere.

3) The complete Hamiltonian $K(I)$ is analytic on the space $\ell^1_{3\/2}$.
Our results show that the function $K(I)-d K(0)(I)=-V(I)$
 is smoother and continuously extends to a  larger   space $\ell^2$. A natural question is if the function
 $$
 K(I)-d K(0)(I)- \tfrac12d^2 K(0)(I,I)=K(I)-P_3+3\|I\|_2^2
 $$ is even smoother and continuously extends to a larger  space,
 etc. We do not know the answer.

4)  By Theorem \ref{T1}, $V(I)$ admits a quadratic upper bound in terms of $P_1$.
The estimate \er{VA2} implies the exponential upper bound for $V$ in terms of $\|I\|_2$.
The  bottle neck of our proof which yields the unpleasant exponential factor in \er{VA2}
is the Bernstein inequality, used in Section 3 to prove Lemma \ref{Te}.
We conjective that, in fact, $V(I)$ is bounded by a polynomials of $\|I\|_2$.

Consider the restriction of the function $V(I)$ to $\wt\ell_+^N$ with any $N\ge 1$.
It is known that the corresponding Hessian is non-degenerate:
\[
\lb{02}
\det \rt\{{\pa^2 V(I)\/\pa I_i\pa I_j}\rt\}_{1\le i,j\le N}\ne 0,\qqq \forall \ I\in \ell_N^+.
\]
This result  was proven in \cite{Kri} with serious omissions,  fixed in \cite{BoK}
(see also Appendix 3.6 in \cite{Ku}). Since $V$ is analytic on $\wt\ell_+^N$, then \er{01}
and \er{02} yield that the Hessian of $V|_{\wt\ell_+^N}$ is a positive $N\ts N$ matrix.
Thus  $V$ is convex on $\ell_N^+$. Since $\wt\ell^\iy=\cup \wt\ell^N$ is dense in $\ell^2$, where
$V$ is continuous, we get

\begin{corollary}
\lb{T4}
The function  $V(I)$ is convex on $\ell_+^2$.
\end{corollary}

{\bf Remark} 5).
 By \er{01} and  remark 2, the function $V$ is strictly convex in some vicinity of the origin in $\ell_+^2$ (note that $\ell_+^2$
 is the only phase-space where $V$ is strictly convex).
 We conjecture that it  is strictly  convex   everywhere in $\ell_+^2$.
\medskip

In difference with $V(I)$, the total   Hamiltonian $K(I)$ is not continuous on $\ell^{2}$ since its linear
part $P_3(I)$ is there an unbounded linear functional.
But $P_3(I)$ contributes to equations \er{hameq} the linear rotations
$$
{\pa b_n\/\pa t}=-(2\pi n)^3  b_{-n},\qqq \qqq {\pa b_{-n}\/\pa t}=(2\pi n)^3b_{n},
\qqq n\ge 1.
$$
So the properties of \er{hameq}  essentially are determined by the component $-V(I)$ of the Hamiltonian
$K(I)$. We also note that since
$P_3(I)$ is a bounded linear functional on the space $\ell_{3\/2}^1\ss \ell^2$, then
the complete Hamiltonian $K(I)=P_3-V$ is concave on $\ell_{3\/2}^1$.
The flow of the KdV equation in the action-angle  variables \eqref{Ac_An} is
$$
(I,\f)\to \rt(I, \f(t)=\f+t K'(I)\rt),  \quad t\in\R,
\qqq K'(I)={\pa K(I)\/\pa I}.
$$
Since the function $K$ is concave and analytic on $\ell_{3\/2}^1$, then the flow-maps are twisting:
$$
\lan \f(t; I_{(2)}, \f_{(1)})-\f(t; I_{(1)}, \f_{(1)}), I_{(2)}-I_{(1)}\ran
=t \lan K'(I_{(2)})-K'(I_{(1)}), I_{(2)}-I_{(1)} \ran\le 0\qquad\forall\,t\ge0.
$$
If  the assertion of Remark 5 above holds true, then L.H.S. is $\le -Ct\|I_{(2)}-I_{(1)}\|_2^2$,
where the positive constant $C$ depends on $I_{(2)},I_{(1)}$.

In the finite-dimensional case convexity (and strict  convexity) of an integrable Hamiltonian significantly
simplifies the study of  long time behavior of actions of solutions for perturbed equations. Similar, we  are certain that results of this work  will help to study perturbations of the KdV equation \eqref{KdV}, especially, those which are Hamiltonian. It is important that as a phase space our results suggest the Hilbert space $\ell^2$, rather than a weighted $\ell^1$-space.


\

\vskip 0.25cm
\section {Momentum, quasimomentum and KdV equation}
\setcounter{equation}{0}

\bigskip

{\bf 2.1. Spectrum of the Hill operator.}
We consider the Hill operator $T$ acting in $L^2(\R)$ and given by
$$
 T=-{d^2\/\,dx^2}+q_0+q(x),
$$
where a 1-periodic potential $q$ (with zero mean-value) belongs to the  Sobolev space $\mH_{\a}, \a\ge -1$ and $q_0$ is a constant (so the potential $q_0+q$ may be a distribution).  Below we recall the results from \cite{K2} on the Hill operator with potentials $q\in \mH_{-1}$.   The spectrum of $T$ is absolutely continuous and consists of intervals
   ({\it spectral bands}) $\gS_n$, separated by {\it gaps} $\g_n$ and is given by
$$
\gS_n=[\l^+_{n-1}, \l^-_n ], \qqq \g_n=(\l^-_n,\l^+_n),\qqq  \where
\qqq
    \l^-_{n-1}\le\l^-_n \le\l^+_{n}, \qq  n\ge
    1.
$$
We choose the constant $q_0=q_0(q)$ in such a way that $\l^+_0=0$. Note that a gap-length $|\g_n|\ge 0$ may be zero.
If the n-th gap degenerates, that is $\g_n=\es$, then the corresponding spectral bands $\gS_{n} $ and $\gS_{n+1}$
merge.
The sequence $0=\l_0^+<\l_1^-\le \ \l_1^+\ <\dots$ form the {\it energy spectrum} of $T$ and is the spectrum of the equation $-y''+(q_0+q)y=\l y$ with the 2-periodic  boundary conditions, i.e. $y(x+2)=y(x), x\in \R$. Here the  equality means that $\l_n^-=\l_n^+$ is a
double eigenvalue. The eigenfunctions, corresponding to $\l_n^{\pm}$, have period 1 when $n$ is even,
and they are antiperiodic, i.e., $y(x+1)=-y(x),\ \  x\in \R,$ when $n$ is odd.
\begin{figure}
\tiny
\unitlength=1.00mm
\special{em:linewidth 0.4pt}
\linethickness{0.4pt}
\begin{picture}(108.67,33.67)
\put(41.00,17.33){\line(1,0){67.67}}
\put(44.33,9.00){\line(0,1){24.67}}
\put(108.33,14.00){\makebox(0,0)[cc]{$\Re\l$}}
\put(41.66,33.67){\makebox(0,0)[cc]{$\Im\l$}}
\put(42.00,14.33){\makebox(0,0)[cc]{$0$}}
\put(44.33,17.33){\linethickness{4.0pt}\line(1,0){11.33}}
\put(66.66,17.33){\linethickness{4.0pt}\line(1,0){11.67}}
\put(82.00,17.33){\linethickness{4.0pt}\line(1,0){12.00}}
\put(95.66,17.33){\linethickness{4.0pt}\line(1,0){11.00}}
\put(46.66,20.00){\makebox(0,0)[cc]{$\l_0^+$}}
\put(56.66,20.33){\makebox(0,0)[cc]{$\l_1^-$}}
\put(68.66,20.33){\makebox(0,0)[cc]{$\l_1^+$}}
\put(78.33,20.33){\makebox(0,0)[cc]{$\l_2^-$}}
\put(84.33,20.33){\makebox(0,0)[cc]{$\l_2^+$}}
\put(93.00,20.33){\makebox(0,0)[cc]{$\l_3^-$}}
\put(98.66,20.33){\makebox(0,0)[cc]{$\l_3^+$}}
\put(106.33,20.33){\makebox(0,0)[cc]{$\l_4^-$}}
\end{picture}
\caption{The spectral domain $\C\sm \cup \gS_n$ and the bands $\gS_n=[\l^+_{n-1},\l^-_n], n\ge 1$}
\lb{sS}
\end{figure}
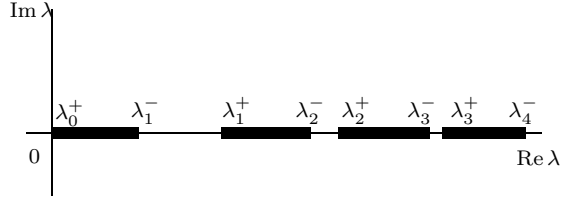

In order to study the actions $I_n, n\ge 1$, we introduce the {\it quasimomentum function}.
We can not introduce the standard fundamental solutions for the operator $T$,
since the perturbation $q$ is too singular if $\a<0$. Instead we use another representation
of $T$. Define a function  $\r(x)$ by
$$
\r(x)=e^{\int_0^xq_*(t)dt},\qqq {\rm where}\qqq q_*\in \mH_0\qqq\qq q_*'=q.
 $$
Consider     the unitary transformation $U:L^2(\R,\r^2 \,dx)\to L^2(\R,dx)$ given by
the multiplication by $\r$. Then $T$ is unitarily equivalent to
$$
T_1y=U^{-1}TUy=-{1\/\r^2}(\r^2y')'+(q_0-q_*^2)y=-y''-2q_*y'+(q_0-q_*^2)y,
\ \ \
$$
acting in $L^2(\R,\r^2 \,dx)$. Note that the norm in this space is equivalent to the original $L^2-$ norm.
This representation clearly is more convenient, since $q_*$ and $q_*^2$ are regular functions.
It is convenient to write the the spectral parameter $\l$ as
$$
\l=z^2.
$$
Let $\vp(x,z)$ and $\vt(x,z)$ be solutions of the equation
\[
\lb{1}
-y''-2q_*y'+(q_0-q_*^2)y=z^2 y, \ \ \ z\in \C ,
\]
satisfying $\vp'(0,z)=\vt(0,z)=1$ and
$\vp(0,z)=\vt'(0,z)=0$. The {\it Lyapunov function} is defined by
\[
\lb{Lf}
\D(z)={1\/2}(\vp'(1,z)+\vt(1,z)).
 \]
 This function is entire and even, i.e.,  $\D(-z)=\D(z)$ for all $z\in \C$.
It is known  that $\D(\sqrt {\l_{n}^{\pm}})=(-1)^n, n\ge 0$ and the function $\D'(z)$ has a unique zero $z_n$ in each gap $[\sqrt {\l_n^-},\sqrt {\l_n^+}]\ss \R_+$ (see e.g., \cite{Kr}, \cite{K7}).

{\bf 2.2. Momentum and quasimomentum.}
Consider a strongly increasing odd sequence $u_n, n\in \Z$,
of real numbers, $u_n=-u_{-n}$,  such that $u_n\to \pm \iy$ as $n\to \pm \iy$, and a non-negative sequence $h=(h_n)_1^\iy\in \ell^\iy$.
We define the following domains
$$
\cK(h)=\C \sm \cup_{n\in\Z }\ol\G_n,\ \ \ \ \ \ \ \  \cK_+(h)=\C  _+\cap \cK(h),
$$
where
$$
\G_0=\es,\qq
\G_n=(u_n -ih_n, u_n+ih_n)=-\G_{-n},\qqq  {\rm and} \qqq \C_+=\{z: \Im z>0\}.
$$
We call $\cK_+(h)$ the "{\it comb}" and denote its points by $k=u+iv$.  Then there exists
a unique conformal mapping $z=z(k)$:
$$
z:\cK_+(h)  \to \C _+ ,\qqq
$$
normalized  by the condition $\ z(0)=0$ and the asymptotics:
\[
\lb{ask1x}
z(iv)=iv+o(v)\qq \as \qq v\to+\iy,\ \ \ {\rm where} \ \ \ z=x+iy,\ \ k=u+iv.
\]
We call $z(k) $ "{\it  the comb mapping}". Define the inverse mapping
\[
\lb{2.k}
k=z^{-1}: {\C _+}\to \cK_+(h),\qqq k(z)=u(z)+iv(z).
\]
This function is continuous in $\C_+$ up to the boundary, i.e., on  the closure $\ol \C_+$. It is convenient to introduce
"gaps" $g_n$, "bands"  $\s_n$ and the "spectrum" $\s$ of the
comb mapping by:
$$
g_n=(z_n^-, z_n^+)=(z(u_n-0), z(u_n+0)),\ \ \ \ \
\s_n=[z_{n-1}^+, z_n^-],\ \ \ \ \ \ \s=\cup\s _n,\qq  g_0=\es, \qq z_0^\pm=0.
$$
Note that the identities $\l_n^\pm={z_n^\pm}^2$ yields
\[
\lb{2.30}
|\g_n|={z_n^+}^2-{z_n^-}^2=|g_n|(z_n^++z_n^-), \qq \forall\  n\ge 1.
\]
Define the momentum domain
$$
\cZ =\C\sm\cup_{n\in \Z} \ol g_n.
$$
 The function $k(z)$ may be  continued from $\C_+$ to the domain $\cZ$ by the symmetry the formula $k(z)=\ol {k}(\ol z),\ \Im z<0$. Thus we obtain a conformal mapping $k: \cZ\to \cK(h)$, called {\it the quasimomentum mapping}
(or shortly the quasimomentum),
which generalizes the classical quasimomentum ( see e.g. \cite{RS}).
A point $z\in \cZ$ is called {\it momentum} and a point $k\in \cK(h)$ is called {\it quasimomentum}.
It is odd, i.e., $k(-z)=-k(z)$, since the domains $\cK(h) $ and $\cZ$  both are invariant under the inversion
$z\to -z$. 

If the spectrum of the comb mapping $k(z)$ has only finite number of open gaps, then  $k(z) $ is called
a {\it finite-gap quasimomentum}. Different properties of the finite-gap quasimomentum (and of more general conformal mappings) were studied by Hilbert one hundred years ago, see in \cite{J}.

The abstract quasimomentum, which we have just defined, is related to the spectral theory
of the Hill operator $T$ by the following construction invented in \cite{MO}.
Namely, let $\{z_n^\pm, n\in\Z\}$ be an odd sequence as above. For $n\ge 0$ denote $\l_n^\pm=(z_n^\pm)^2$.
Then $\{\l_n^\pm, n\ge 0\}$ is the energy spectrum of the Hill operator $T$ with a potential
$q_0+q$, where $q\in \mH_\a, \a\ge 0$, if and only if the corresponding comb domain $\cK(h)$ is such that
$u_n=\pi n, n\in \Z$ and $h=(h_n)_1^\iy\in \ell_{\a+1}^2$. Moreover, in this case $\cos k(z)= \D(z)$ is the Lyapunov function for $T$.

In \cite{K2} the construction  was  generalized for  potentials from $\mH_{-1}$, see below Theorem \ref{T1.1}.

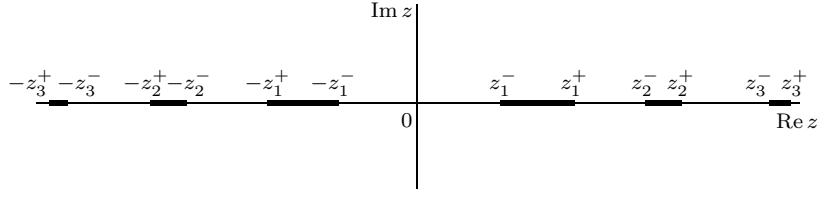
\begin{figure}
\tiny
\unitlength=1mm
\special{em:linewidth 0.4pt}
\linethickness{0.4pt}
\begin{picture}(120.67,34.33)
\put(20.33,21.33){\line(1,0){100.33}}
\put(70.33,10.00){\line(0,1){24.33}}
\put(69.00,19.00){\makebox(0,0)[cc]{$0$}}
\put(120.33,19.00){\makebox(0,0)[cc]{$\Re z$}}
\put(67.00,33.67){\makebox(0,0)[cc]{$\Im z$}}
\put(81.33,21.33){\linethickness{2.0pt}\line(1,0){9.67}}
\put(100.33,21.33){\linethickness{2.0pt}\line(1,0){4.67}}
\put(116.67,21.33){\linethickness{2.0pt}\line(1,0){2.67}}
\put(60.00,21.33){\linethickness{2.0pt}\line(-1,0){9.33}}
\put(40.00,21.33){\linethickness{2.0pt}\line(-1,0){4.67}}
\put(24.33,21.33){\linethickness{2.0pt}\line(-1,0){2.33}}
\put(81.67,24.00){\makebox(0,0)[cc]{$z_1^-$}}
\put(91.00,24.00){\makebox(0,0)[cc]{$z_1^+$}}
\put(100.33,24.00){\makebox(0,0)[cc]{$z_2^-$}}
\put(105.00,24.00){\makebox(0,0)[cc]{$z_2^+$}}
\put(115.33,24.00){\makebox(0,0)[cc]{$z_3^-$}}
\put(120.00,24.00){\makebox(0,0)[cc]{$z_3^+$}}
\put(59.33,24.00){\makebox(0,0)[cc]{$-z_1^-$}}
\put(50.67,24.00){\makebox(0,0)[cc]{$-z_1^+$}}
\put(40.33,24.00){\makebox(0,0)[cc]{$-z_2^-$}}
\put(34.67,24.00){\makebox(0,0)[cc]{$-z_2^+$}}
\put(26.00,24.00){\makebox(0,0)[cc]{$-z_3^-$}}
\put(19.50,24.00){\makebox(0,0)[cc]{$-z_3^+$}}
\end{picture}
\caption{$z$-domain $\cZ=\C\sm\cup g_n$, where $z=\sqrt{\l}$ and momentum gaps $g_n=(z_n^-,z_n^+)$}
\lb{z}
\end{figure}
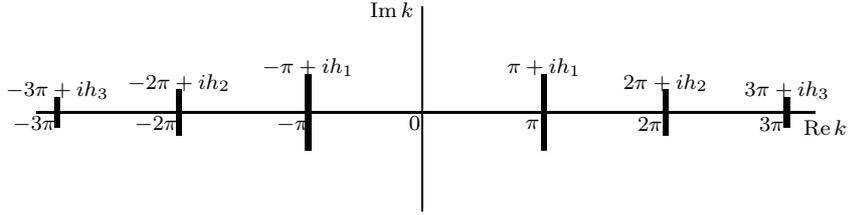
\begin{figure}
\tiny
\unitlength=1mm
\special{em:linewidth 0.4pt}
\linethickness{0.4pt}
\begin{picture}(120.67,34.33)
\put(20.33,20.00){\line(1,0){102.33}}
\put(71.00,7.00){\line(0,1){27.00}}
\put(70.00,18.67){\makebox(0,0)[cc]{$0$}}
\put(124.00,18.00){\makebox(0,0)[cc]{$\Re k$}}
\put(67.00,33.67){\makebox(0,0)[cc]{$\Im k$}}
\put(87.00,15.00){\linethickness{2.0pt}\line(0,1){10.}}
\put(103.00,17.00){\linethickness{2.0pt}\line(0,1){6.}}
\put(119.00,18.00){\linethickness{2.0pt}\line(0,1){4.}}
\put(56.00,15.00){\linethickness{2.0pt}\line(0,1){10.}}
\put(39.00,17.00){\linethickness{2.0pt}\line(0,1){6.}}
\put(23.00,18.00){\linethickness{2.0pt}\line(0,1){4.}}
\put(85.50,18.50){\makebox(0,0)[cc]{$\pi$}}
\put(54.00,18.50){\makebox(0,0)[cc]{$-\pi$}}
\put(101.00,18.50){\makebox(0,0)[cc]{$2\pi$}}
\put(36.00,18.50){\makebox(0,0)[cc]{$-2\pi$}}
\put(117.00,18.50){\makebox(0,0)[cc]{$3\pi$}}
\put(20.00,18.50){\makebox(0,0)[cc]{$-3\pi$}}
\put(87.00,26.00){\makebox(0,0)[cc]{$\pi+ih_1$}}
\put(56.00,26.00){\makebox(0,0)[cc]{$-\pi+ih_1$}}
\put(103.00,24.00){\makebox(0,0)[cc]{$2\pi+ih_2$}}
\put(39.00,24.00){\makebox(0,0)[cc]{$-2\pi+ih_2$}}
\put(119.00,23.00){\makebox(0,0)[cc]{$3\pi+ih_3$}}
\put(23.00,23.00){\makebox(0,0)[cc]{$-3\pi+ih_3$}}
\end{picture}
\caption{$k$-plane and cuts $\G_n=(\pi n-ih_n,\pi n+ih_n), n\in\Z$}
\lb{k}
\end{figure}

\begin{figure}
\unitlength 1mm 
\linethickness{0.4pt}
\ifx\plotpoint\undefined\newsavebox{\plotpoint}\fi 
\begin{picture}(119.75,66.5)(0,0)
\put(15,17.25){\line(1,0){104.75}}
\qbezier(40.5,17.25)(40.5,30.375)(52.25,35.75)
\qbezier(52.25,35.75)(58.25,38.375)(64.25,35.75)
\qbezier(64.25,35.75)(76.17,30.375)(76,17.25)
\qbezier(40.5,17.25)(40.5,45.125)(58.25,46.5)
\qbezier(76,17.25)(76.17,45.125)(58.25,46.5)
\put(38.75,12.00){\makebox(0,0)[cc]{$z_n^-$}}
\put(76,12){\makebox(0,0)[cc]{$z_n^+$}}
\put(65.25,38.5){\makebox(0,0)[cc]{$v_n$}}
\put(70,44.25){\makebox(0,0)[cc]{$v$}}
\put(18,17.25){\line(0,-1){1.00}}
\put(95,17.25){\line(0,-1){1.00}}
\put(57,17.25){\line(0,-1){1.00}}
\put(18,12.00){\makebox(0,0)[cc]{$z_{n-1}^+$}}
\put(95,12.00){\makebox(0,0)[cc]{$z_{n+1}^-$}}
\put(57,12.00){\makebox(0,0)[cc]{$z_{n}$}}
\end{picture}
\caption{The graph of $v(z+i0), \ z\in g_n\cup \s_n\cup \s_{n+1}$
and $|h_n|=v(z_n+i0)>0$}
\lb{grafv}
\end{figure}
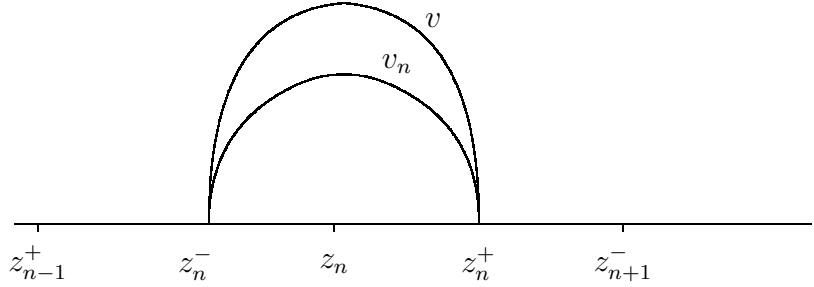

Despite the objects, treated by Theorem \ref{T1}, are defined in terms of Hill operators with
periodic potentials, for the proofs in Sections 3-4 below we need the quasimomentum mapping
$k(z), k=u+iv, z=x+iy$, corresponding
to  general odd sequences $\{u_n\}$. Now we summarize their basic properties, refering for a proof to \cite{KK2} \cite{K1}, \cite{K2}, \cite{K4}, \cite{L1}-\cite{L3}, \cite{MO}, \cite{M}.

\no 1) {\it $v(z)\ge \Im z>0$ and $v(z)=-v(\ol z)$ for all $z\in \C_+$ and
\[
k(-z)=-k(z), \qqq k(z)=\ol k(\ol z) ,\qq \all \ z\in \cZ,
\]

\no 2) $v(z)=0$ for all $z\in \s_n=[z_{n-1}^+,z_n^-], n\in \Z$.

\no 3)  If some $g_n\ne \es, n\in \Z$,  then
\[
\lb{prqx}
h_n\ge v(z+i0)=-v(z-i0)>0, \qqq v''(z+i0)<0 \qqq  \forall \  z\in g_n,
\]
see Fig. 4. The function $v(z+i0)|_{g_n}>0$  attains its maximum at a point $z_n\in g_n$, where
\[
\lb{2.04}
h_n=v(z_n+i0), \qq v'(z_n)=0.
\]
Moreover,
\[
\lb{2.40}
v=0 \qqq {\rm on}\qq \R\sm \cup_{n\in \Z}g_n,
\]
\[
\lb{prq}
v(z+i0)>v_n(z)=|(z-z_n^-)(z-z_n^+)|^{1\/2}>0, \qqq \qqq  \forall \  z\in g_n,
\]
\[
\lb{em-1}
|g_n|\le 2h_n,\qqq |\s_n|\le u_n-u_{n-1},\qqq  \forall \  n\in \Z.
\]
\no 4) $u'(z)>0$ on  each $\s_n$, and
\[
\lb{2.77}
u(z)=\pi n, \qqq \forall \ z\in g_n\ne \es, n\in \Z.
\]

\no 5) The function $k(z)$ maps a horizontal cut (a gap ) $\ol g_n$  onto the vertical cut $\ol \G_n$  and maps
a spectral band $\s_n$ onto the segment $[\pi (n-1), \pi n]$, for all $n\in\Z$.

\no 6) The following asymptotics hold true:
\[
\lb{aen}
z_n^\pm=\pi n+ o(1) \qqq \as \qq n\to \iy.
\]
\no  7) If $h\in \ell^2$ and $\inf_{n\ge 1} (u_{n+1}-u_n)>0$, then $v(z+i0), z\in \R$ belongs to $L^1(\R)$ and the following identity holds true:
\[
\lb{Kq}
k(z)=z+{1\/\pi}\int_{\bigcup_{n\in \Z} g_n} {v(t)\/t-z}dt,\qqq \forall z\in \cZ.
\]
}


For additional properties of the comb mapping  $z(k)$, see \cite{K1}-\cite{K6}, \cite{L1}-\cite{L3}.

{\bf 2.3. Quasimomentum and the KdV Hamiltonian.}
Recall that we choose the constant $q_0\ge 0$ in such a way that  $\l_0^+=0$.
If $q\in \mH_0(\T)$, then the quasimomentum $k(\cdot)$ has asymptotics
\[
\lb{ak}
k(z)=z-{Q_0\/z}-{Q_2+o(1)\/z^3}\qqq as \qq \Im z\to \iy,
\]
see \cite{K2}.
If $q,q'\in \mH_0$, then the asymptotics \er{ak} may be improved:
\[
\lb{Tm-1}
k(z)=z-{Q_0\/z}- {Q_2\/z^3}-{Q_4+o(1)\/z^5} \qq as \qq z\to +i\iy,
\]
and
\[
\lb{ask1}
k^2(z)=\l-S_{-1}-{S_0 \/\l} -{S_1+o(1) \/ \l^2} \qqq
\as \ \ \l=z^2, \qq z\to +i\iy,
\]
 where
\[
\lb{2.130}
 Q_j={1\/\pi}\int_\R z^jv(z+i0)\,dz\ge 0,\qq j\ge 0,
\qqq S_j={4\/\pi}\int_{0}^\iy z^{2j+1}u(z)v(z+i0)\,dz, \qq j\ge -1.
\]
Note that $Q_{2j+1}=0$ for all $j\ge 0$ by the symmetry.  The involved quantities $Q_j$ and $S_j$ are defined by  converging integrals  (see \cite{KK2}, \cite{K5}), and satisfy  the following identities
\[
\lb{i1}
q_0(q)=S_{-1}=2Q_0 \qqq {\rm if }\qqq q\in \mH_{-1},
\]
\[
\lb{H1}
H_1(q)=\int_0^1q^2(x)\,dx=2P_1=4S_0=8Q_2-4Q_0^2 \qqq {\rm if }\qqq q\in \mH_{0},
\]
\[
\lb{H2}
 H_2(q)=8(S_1-S_{-1}S_0),\qqq S_{1}+2Q_0Q_2=2Q_4 \qqq {\rm if }\qqq q\in \mH_{1},
\]
\[
\lb{ask2}
8Q_2=\|q\|^2+q_0^2,\qqq  2^4Q_4=H_2(q+q_0).
\]
See  \cite{K2}, \cite{K5}.

{\bf 2.4. The KdV actions.}
The components $I_n$ of the action vector $I=(I_n)_1^\iy$ (see \er{1.22})
may be calculated with the help of a general formula due to Arnold,
which in the KdV-case takes the form
\[
I_n={(-1)^{n+1}2\/\pi}\int_{g_n}{z^2\D'(z)\,dz\/|\D^2(z)-1|^{1\/2}}\ge
0,\qqq n\ge 1,
\]
see \cite{FM}. These integrals may be re-written, using the quasimomentum.
Indeed, since $\sin k(z)=\sqrt{1-\D^2(z)}$, then
$$
I_n=-{1\/\pi i}\int_{c_n}z^2{\D'(z)\/\sin k(z)}\,dz,
$$
where $c_n$ is a contour around $g_n$. It is convenient to introduce  contours $\c_n$ around $\G_n$ by
\[
\lb{decc} \c_n=\rt\{k\in \cK(h): {\rm dist}\ \ (k, \G_n)={\pi\/4}\rt\}\ss
\cK(h),\qqq n\ge 1,
\]
and define the contours $c_n$  as
$$
c_n=z(\c_n)\ss \cZ, \qqq \forall \ n\ge 1.
$$
 The differentiation of $\D(z)=\cos k(z)$ gives $k'(z)=-\D'(z) / {\sin k(z)}$. This yields
\[
\lb{deac}
I_n={1\/i\pi}\int_{c_n}z^2k'(z)\,dz=-{2\/i\pi}\int_{c_n}zk(z)\,dz={4\/\pi}\int_{g_n}zv(z+i0)\,dz\ge 0,
\]
since   on $g_n$ the function $k=\pi n+iv$  and $v$ satisfies \er{prqx}.
This representation for $I_n$ is convenient and is crucial for our work. In particular, below in Lemma
 \ref{Te} we derive from \er{deac} the following two-sided estimates:
 $$
 {2\/3\pi}h_n|\g_n|<I_n\le {2h_n|\g_n|\/\pi },\qqq {\rm if} \qq |\g_n|>0.
 $$
Using \er{deac}  jointly with \er{2.77} and \er{2.40}, we easily see that
\[
\lb{2.50}
P_3=\sum_{n\ge 1}(2\pi n)^3I_n={32\/\pi}\int_0^\iy zu^3(z)v(z+i0)\,dz.
\]
Recall that $P_3(I)$ is the linear in $I$ part of the Hamiltonian $H_2$, see \er{HVid}.

{\bf 2.5. Marchenko-Ostrovski construction for potentials $q\in \mH_{-1}$.}
The Marchenko-Ostrovski construction, described in Section 2.1, defines the mapping $q\to h$, acting from
 $\mH_{j-1}$ into $\ell^2_j, j=0,1,....$. The results below are proven in \cite{MO} for
 $j\ge 1$ and in \cite{K2} for $j=0$.

\begin{theorem}
\lb{T1.1}
The mapping $q\to h$ acting from $\mH_{j-1}$ into $\ell^2_j, j=0,1$ is a surjection.
It satisfies the following estimates
\[
\lb{1.6}
\|q\|_{-1}\le 2\|h\|_2(1+4\|h\|_2),\ \ \ \ \|h\|_2\le 3\|q\|_{-1}(1+2\|q\|_{-1})^2,
\ \ \ \forall\  q\in \mH_{-1},
\]
where $\|q\|_{-1}=\|q\|_{\mH_{-1}}$.
For each $h\in\ell^2$ there exists a function $q\in \mH_{-1}, $
and a unique conformal mapping $k(\cdot,h):\cZ\to K(h)$ defined in \er{2.k}. Moreover,
\[
\lb{1.16}
\cos k(z,h)=\D(z,h),\ \ \ \ \ z\in\cZ,
\]
where $\D(z,h)$ is the Lyapunov function for $q$, and $k(z)$ satisfy
\[
\lb{1.17}
k(z,h)=z-{Q_0+o(1) \/ z}\ \ \ \ \ \ {\rm as}\ \ \ \ \ z\to i\iy,
\]
\[
\lb{1.19}
k(z_n^\pm,h)=\pi n\pm i0,\qqq  k(z_n\pm i0,h)=\pi n\pm ih_n,\ \ \  n\ge 1.
\]
In particular, the real numbers ${z_n^\pm}^2$, satisfying \er{1.19}, form the energy spectrum of the operator $T$.
Furthermore,  if a sequence $h^\n, \n\ge 1$ converges strongly in $\ell^2$ to $h$
as $\n\to\iy$, then $\D(z,h^\n)\to\D(z,h)$ uniformly on bounded subsets of $\C$.
\end{theorem}

In order to prove our main result, Theorem  \ref{T1},  we use Theorem  \ref{T1.1} to reformulate it
as questions from the conformal mapping theory in terms of quasimomentum of the Hill operator.
To proceed we need auxillary results from previous work of the first author:

\begin{lemma}
\lb{Tae}
Let $h\in \ell^2$ and let each $u_n=\pi n, n\ge 1$. Then following estimates hold true:
\[
\lb{11} {\pi \/4}Q_0\le\|h\|_2^2 \le
 {\pi^2\/ 2}\lt(1+{\sqrt 2\/\pi}Q_0^{1\/2}\rt)Q_0,
\]
\[
\lb{12}
 \|\r\|_2^2\le (16)^2 Q_0,
\]
\[
\lb{Te1}
{\|h\|_{\iy}^2\/2}\le Q_0\le {2\/\pi}\int_0^\iy {zv(z)\,dz\/u(z)}=\sum_{j\ge 1}{I_j\/2\pi j}=P_{-1},
\]
\[
\lb{Teg} \|h\|_{\iy}\le {4\/\pi}\sum_{j=1}^\iy{|\g_j|\/2\pi j}.
\]

\end{lemma}
\no {\bf Proof.}
Estimates \er{11}, \er{12} were proved in Theorem 2.1 from \cite{K1}.
The first estimate in \er{Te1} $\|h\|_\iy^2\le 2Q_0$ was established in
\cite{K6}. The second estimate in \er{Te1} $Q_0\le {2\/\pi}\int_0^\iy {zv(z)\,dz/u(z)}$
was proved in \cite{K5} (see p. 398). The identities
\er{deac},  \er{2.40} and \er{2.77} imply
$$
{2\/\pi}\int_0^\iy
{zv(z)\,dz\/u(z)}=\sum_{j\ge 1}{I_j\/2\pi j},
$$
and we get \er{Te1}. Finally, using  \er{Te1}, \er{2.04} and  \er{2.77} we obtain
\[
\lb{Te11} \|h\|_\iy^2\le {4\/\pi}\int_0^\iy {zv(z)\,dz\/u(z)}\le
{4\/\pi}\sum_{j\ge 1} h_j\int_{g_j}{z\,dz\/u(z)}= {2\/\pi}\sum_{j\ge 1} h_j{|\g_j|\/\pi
j}\le {2\|h\|_\iy\/\pi}\sum_{j\ge 1}{|\g_j|\/\pi j},
\]
which gives \er{Teg}.
\BBox

\bigskip

\vskip 0.25cm
\section {Local estimates }
\setcounter{equation}{0}

\bigskip

In this section we derive estimates for $h_n, I_n$ and $|\g_n|$ with a fixed $n\ge 1$.
We use the following constants
\[
\lb{3.0}
C_-=e^{\sqrt{P_{-1}}},\qqq C_I=1+\sqrt{P_{-1}}\qqq C_0=\ch \|h\|_\iy\le e^{\sqrt{2P_{-1}}},
\]
where the inequality follows from lemma below.

\begin{lemma}
\lb{Te}
Let $h\in \ell^2$ and let each $u_n=\pi n, n\ge 1$. Then for each $n\ge 1$ the following estimates hold true:
\[
\lb{Te3}
{2\/3\pi}h_n|\g_n|<{2\/3\pi}h_n|g_n|(z_n +z_n^-+z_n^+)<
I_n\le {2h_n|\g_n|\/\pi },\qqq {\rm if} \qq |\g_n|>0,
\]
\[
\lb{Te41}
z_n^\pm \le \pi n+\sum_{j= 1}^n|g_j|,
\]
\[
\lb{Te42}
\pi n\le 2z_n^\pm +{\|\r\|_2^2\/\pi}, \qqq \qqq \r=(\r_n)_1^\iy,\qq \r_n=\pi-|\s_n|,
\]
\[
\lb{Te4}
2n\le C_0z_n^\pm, \qqq  \qqq h_n\le {\sqrt{C_0}\/2}|g_n|,
\]
\[
\lb{Te5} 2\pi n h_n^2\le \sqrt{C_0}{3\pi\/2}I_n+2{\|\r\|^2\/\pi}h_n^2,
\]
\[
\lb{Te6} {1\/C_I}{|\g_n|\/4\pi n}\le h_n\le {\pi
C_0^{3\/2}|\g_n|\/8\pi n}.
\]
\[
\lb{gIn} {1\/3\pi C_I}{|\g_n|^2\/(2\pi n)}\le I_n\le {C_0^{3\/2}\/2}{|\g_n|^2\/(2\pi n)},\qqq
\]
\[
\lb{hIn} {8C_0^{-{3\/2}}\/3\pi^2}(2\pi n)|h_n|^2\le I_n\le  8nC_Ih_n^2.
\]
\end{lemma}
\no {\bf Proof.}
We show \er{Te3}. Using  \er{prqx}, \er{2.04} and  standard convexity arguments (see Fig. 5) we have
\[
\begin{aligned}
\lb{f-}
v(z_n^-+t+i0)\ge f_-(t)=t{h_n\/\ve_-}, \qqq t\in (0,\ve_-),\qq \ve_-=z_n-z_n^->0.
\end{aligned}
\]
\begin{figure}
\tiny
\unitlength 1mm 
\linethickness{0.4pt}
\ifx\plotpoint\undefined\newsavebox{\plotpoint}\fi 
\begin{picture}(126.775,106.662)(0,0)
\put(3.75,24.85){\line(1,0){123.025}}
\qbezier(25.275,24.675)(59.575,106.662)(93.175,24.5)
\multiput(93.175,24.5)(-.03372855701,.04114530777){991}{\line(0,1){.04114530777}}
\multiput(25.625,25.025)(.033726647,.03992133727){1017}{\line(0,1){.03992133727}}
\put(59.575,65.275){\line(0,-1){40.25}}
\put(24.925,21.5){\makebox(0,0)[cc]{$z_n^-$}}
\put(92.65,21.5){\makebox(0,0)[cc]{$z_n^+$}}
\put(59.575,21){\makebox(0,0)[cc]{$z_n$}}
\put(41.725,26.3){\makebox(0,0)[cc]{$\ve_-$}}
\put(73.925,26.3){\makebox(0,0)[cc]{$\ve_+$}}
\put(45.575,44.45){\makebox(0,0)[cc]{$f_-$}}
\put(71.65,44.275){\makebox(0,0)[cc]{$f_+$}}
\put(63.075,68.25){\makebox(0,0)[cc]{$v(z_n+i0)=h_n$}}
\put(85.3,53.2){\makebox(0,0)[cc]{$v(z+i0)$}}
\end{picture}
\caption{\footnotesize The graphs of $v(z+i0)$ and $f_\pm$}
\lb{Fi5}
\end{figure}
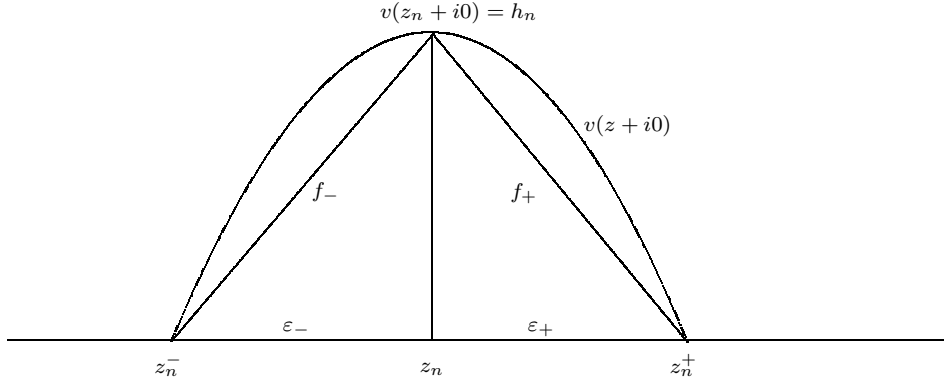
 This yields
$$
I_n^-={4\/\pi}\int_{z_n^-}^{z_n} z v(z)\,dz\ge {4\/\pi}\int_0^{\ve_-}(z_n^-+t)f_-(t)dt=
{4\/\pi}{h_n\/\ve_-}\rt(z_n^-{\ve_-^2\/2}+ {\ve_-^3\/3}\rt)
={2\/\pi}h_n\ve_-\rt(z_n-{\ve_-\/3}\rt).
$$
Let $f_+(t)=(\ve_+-t){h_n\/\ve_+}, t\in (0,\ve_+), \ve_+=z_n^+-z_n>0$.
A similar argument gives
$$
I_n^+={4\/\pi}\int_{z_n}^{z_n^+} z v(z)\,dz\ge {4\/\pi}\int_0^{\ve_+}(z_n^++t-\ve_+)f_+(t)dt=
{4\/\pi}{h_n\/\ve_+}\rt(z_n^+{\ve_+^2\/2}-{\ve_+^3\/3}\rt)
={2\/\pi}h_n\ve_+\rt(z_n+{\ve_+\/3}\rt).
$$
Denoting $z_n^0={1\/2}(z_n^++z_n^-)$, we obtain
$$
I_n=I_n^-+I_n^+>
{2\/\pi}h_n\rt(\ve_-(z_n-{\ve_-\/3}) +\ve_+(z_n+{\ve_+\/3})\rt)=
{2\/\pi}h_n\rt(z_n |g_n|+{\ve_+^2-\ve_-^2\/3}\rt)
$$
$$
={2\/\pi}h_n|g_n|\rt(z_n +2{z_n^0-z_n\/3}\rt)={2\/3\pi}h_n|g_n|(z_n +2z_n^0)\ge {2\/3\pi}h_n|g_n|2z_n^0={2\/3\pi}h_n|\g_n|.
$$
This implies the first two estimates in \er{Te3}.
Using \er{prqx} for all $z\in g_n$, we get $I_n<{4\/\pi}h_n\int_{g_n} z\,dz={2\/\pi }h_n|\g_n|$,
which gives us the last estimate in \er{Te3}.

We show \er{Te41}. It is clear that
\[
\lb{123}
z_n^+=\pi n+\sum_1^n \rt(|g_j|-(\pi -|\s_j|)\rt)
, \qqq \all \ n\ge 1.
\]
Since by \er{em-1}, $\r_j=\pi-|\s_j|\ge 0$  and $|g_n|=z_n^+-z_n^-$,  we get \er{Te41}.

We show \er{Te42}. Using identities \er{123} we obtain
$$
\pi n\le z_n^\pm +\sum_1^n \r_j\le z_n^\pm +n^{1\/2}\|\r\|_2\le z_n^\pm +{\pi n\/2}+{\|\r\|_2^2\/2\pi},
$$
which yields \er{Te42}.

In order to prove \er{Te4} we use an argument from \cite{MO}. This  is a weak point in our proof, which
gives the exponential factor in \er{Te4} and later in \er{VA2}.
The Taylor formula implies
\[
\lb{Tei1}
2=|\D(z_n^-)-\D(z_{n-1}^+)|\le |\D'(\wt z_n)||\s_n|,
\]
for some $\wt z_n\in \s_n=[z_{n-1}^+, z_{n}^-]$ and all $n\ge 1$.
Using the Bernstein inequality for the bounded exponential type functions (see e.g., \cite{S}) we obtain
\[
\lb{Be1} \sup_{z\in \R} |\D'(z)|\le \sup_{z\in \R} |\D(z)|=C_0
=\ch \|h\|_\iy\le  e^{\|h\|_\iy}.
\]
Combining \er{Tei1} and \er{Be1} we get $2n\le C_0z_n^\pm$ for all $n\ge 1$, which yields the first estimate in \er{Te42}.

Let $n$ be even and let $|z_n^--z_n|\le |g_n|/2$ (for other cases the proof is similar). The identity $\ch {h_n}=\D(z_n)$ (which follows from \er{1.19}) and the  Taylor formula imply
\[
\lb{Te2}
{h_n^2\/2}\le \ch {h_n}-1=\D(z_n)-1={1\/2}\D''(\wt z_n^-)(z_n^--z_n)^2
\]
 for some $\wt z_n^-\in (z_n^-,z_n)$.  Using again the Bernstein inequality, we obtain
\[
\lb{Be2}
\sup_{z\in \R} |\D''(z)|\le \sup_{z\in \R} |\D(z)|=C_0.
 \]
 Then combining \er{Te2} and \er{Be2}
 we get $h_n^2\le {C_0\/4}|g_n|^2$, which gives  the second estimate in \er{Te4}.

We show  \er{Te5}. Using \er{Te42}, \er{Te4}, \er{Te3}, \er{2.30} we obtain
$$
\pi n h_n^2\le (z_n^-+z_n^+)h_n^2 +{\|\r\|_2^2\/\pi}h_n^2 \le (z_n^-+z_n^+){\sqrt{C_0}\/2}|g_n|h_n+{\|\r\|_2^2\/\pi}h_n^2
$$$$
={\sqrt{C_0}\/2}|\g_n|h_n+{\|\r\|_2^2\/\pi}h_n^2\le {\sqrt{C_0}\/2}{3\pi\/2}I_n+{\|\r\|_2^2\/\pi}h_n^2,
$$
which yields \er{Te5}.

We show  \er{Te6}.
Using \er{2.30}, \er{Te41}, \er{Te1} and \er{em-1} we obtain
$$
{|\g_n|\/4\pi n} ={(z_n^-+z_n^+)|g_n|\/4\pi n} \le \rt(1+{\|g\|_\iy\/ \pi}\rt)
h_n \le \rt(1+{\sqrt{8P_{-1}}\/ \pi}\rt) h_n\le C_I h_n.
$$
Recalling that $C_0=\ch \|h\|_\iy$ and using \er{Te4} and \er{2.30}, we obtain
$$
h_n\le {C_0^{1\/2}|g_n|\/2}=
 {C_0^{1\/2}|\g_n|\/2(z_n^-+z_n^+)}\le {\pi
 C_0^{3\/2}|\g_n|\/8\pi n},
$$
and \er{Te6} is proven.

Estimates \er{Te6} and \er{Te3} imply the first estimate in \er{gIn}:
$$
{|\g_n|^2\/(2\pi n) }\le 2C_Ih_n{|\g_n|}\le 3\pi C_I  I_n.
$$
Combining the last estimate in \er{Te6} and \er{Te3} we obtain the second estimate in \er{gIn}.

We show  \er{hIn}. Using \er{Te3}  and  \er{Te6}  we obtain
$$
I_n\le {2h_n|\g_n|\/\pi }\le 8nC_Ih_n^2,
$$
which yields the second estimate in \er{hIn}. Using \er{Te6} and \er{Te3}, we obtain
$$
2\pi nh_n^2\le {\pi\/4} C_0^{3\/2}h_n|\g_n|\le  {3\pi^2\/8} C_0^{3\/2}I_n,
$$
and get the first.
 \BBox

For any $h\in \ell^\iy$ we define  integrals  $V_n$ as
\[
\lb{dVn} V_n={8\/\pi}\int_{g_n} zv^3(z)\,dz\ge 0,\qqq n\ge 1.
\]
These quantities are important for our argument since, as we show below,
 $V=\sum_{n\ge 1}(4\pi n)V_n$ for $I\in \ell^2$.

\begin{lemma}
\lb{TgapI}
Let $h\in \ell^\iy$ and let each $u_n=\pi n, n\ge 1$. Then for $n\ge 1$ the following relations hold true:
\[
\lb{VhVn}
{1\/5}h_n^2I_n\le {2\/5\pi}h_n^3|\g_n|\le
{2\/5\pi}h_n^3|g_n|(3z_n +2z_n^0)\le V_n\le 2h_n^2I_n,
\]
\[
\lb{VnI} {4\/\pi i}\int_{c_n}zk^4(z)\,dz=(4\pi n)V_n-(2\pi n)^3I_n.
\]

\end{lemma}
\no {\bf Proof}.
We show \er{VhVn}.  Let $g_n\ne \es$.  Using \er{f-} we get $v(z_n^-+t+i0)\ge f_-(t)=  t h_n / {\ve_-}$,  $t\in (0,\ve_-)$, where $\ve_-=z_n-z_n^->0$. Therefore
$$
V_n^-:={8\/\pi}\int_{z_n^-}^{z_n} z v^3(z)\,dz\ge {8\/\pi}\int_0^{\ve_-}(z_n^-+t)f_-^3(t)dt=
{8\/\pi}{h_n^3\/\ve_-^3}\rt(z_n^-{\ve_-^4\/4}+ {\ve_-^5\/5}\rt)
$$$$
={2\/\pi}h_n^3\ve_-\rt(z_n^-+ {4\ve_-\/5}\rt)={2\/\pi}h_n^3\ve_-\rt(z_n-{\ve_-\/5}\rt).
$$
Similar argument yields $v(z_n+t+i0)\ge f_+(t)=(\ve_+-t) h_n /{\ve_+}$,  $t\in (0,\ve_+)$, where $\ve_+=z_n^+-z_n>0$.
Thus
$$
V_n^+={8\/\pi}\int_{z_n}^{z_n^+} z v^3(z)\,dz\ge {8\/\pi}\int_0^{\ve_+}(z_n^++t-\ve_+)f_+^3(t)dt=
{8\/\pi}{h_n^3\/\ve_+^3}\rt(z_n^+{\ve_+^4\/4}-{\ve_+^5\/5}\rt)
$$$$
={2\/\pi}h_n^3\ve_+\rt(z_n^+-{4\ve_+\/5}\rt)={2\/\pi}h_n^3\ve_+\rt(z_n+{\ve_+\/5}\rt).
$$
Summing these relations  we obtain
$$
V_n=V_n^-+V_n^+>
{2\/\pi}h_n^3\rt(\ve_-(z_n-{\ve_-\/5}) +\ve_+(z_n+{\ve_+\/5})\rt)=
{2\/\pi}h_n^3\rt(z_n |g_n|+{\ve_+^2-\ve_-^2\/5}\rt)
$$
$$
={2\/\pi}h_n^3|g_n|\rt(z_n +2{z_n^0-z_n\/5}\rt)={2\/5\pi}h_n^3|g_n|(3z_n +2z_n^0).
$$
Using \er{prqx}, we get $V_n\le (8h_n^2/\pi)\int_{g_n} zv(z,h)\,dz=2h_n^2I_n$, which yields the last two
estimates in \er{VhVn}. The second follows from \er{2.30} and the first follows from the last estimate in \er{Te3}.

 Using \er{prqx}, \er{2.77} and \er{deac} we obtain
$$
{4\/\pi i}\int_{c_n}zk^4(z)\,dz={4\/\pi i}\int_{c_n}z(u^2-v^2+2iuv)^2\,dz=-{8\/\pi i}\int_{g_n}2z(u^2-v^2)(2iuv)\,dz
$$
$$
={32\/\pi }\int_{g_n} z(v^2-u^2)uv\,dz=(4\pi n)V_n-(2\pi n)^3I_n.
$$
This proves \er{VnI}.
 \BBox
 \medskip

{\no \bf Remark.} 1) In
particular, \er{gIn}  yields that $I\in \ell^2$ iff $\g\in \ell_{-1}^4$.

2)  Due to \cite{K2}, for any $N_0>1$ and $\ve \in (0,{1\/4})$ there exists a potential $q\in \mH_{-1}$
 such that $|\g_n|=n^\ve, $ for all $n>N_0$. Then
$\g=(|\g_n|)_1^\iy\in \ell_{-1}^4$  and \er{gIn}
gives that $I\in \ell^2$. It is clear that $q\notin L^2(\T)$,
since the gap length $|\g_n|$ is increasing.

Note that if $I\in \ell^2$, then \er{Te6} gives $\sum_{n\ge 1}|\g_n|^4/n^2<\iy$, which yields $(|\g_n|)_1^\iy\in
\ell_{-{1\/2}}^2$. This and the standard relationship between the gap lengths and the Fourier coefficients
of potential imply $q\in \mH_{-{1\/2}}$ (e.g., see \cite{K2}, where an analogy of this relation is established for
$(|\g_n|)_1^\iy\in \ell_{-{1}}^2$  and for $(|\g_n|)_1^\iy\in \ell^2$ ).

3) Relations \er{gIn} show that asymptotically the actions $I_n$ are equivalent to the weighted squared gap-length
${|\g_n|^2/ 2\pi n}$. It is known  that the gap-length $|\g_n|$ is asymptotically  equivalent
to the module of the Fourier coefficients $|\hat q_n|$ of the potential $q$ (see \cite{MO}, \cite{M}
for the equivalence and see \cite{K4} for the corresponding estimates).
The asymptotical  equivalences \ $
nI_n\sim |\g_n|^2\sim |\hat q_n|^2
$ \
and the corresponding estimates are important for the spectral theory of the Hill operator $T$
and the theory of KdV.

\bigskip

\vskip 0.25cm
\section {Proof of Theorem \ref{T1}}
\setcounter{equation}{0}

\bigskip


We remind that $V(I)=P_3-H_2$ is the non-linear part of the KdV Hamiltonian, written as a function of actions,
see \er{1.70} and  \er{HVid}. Identities \er{H2} and \er{deac}  represent  $H_2$ and $I_j$ as integrals in terms of the
quasimomentum $k=u+iv$. They allow to write $V$ in a similar form. We start with an integral representation
for $V$ for the case of smooth potentials.

\begin{lemma}
\lb{TQ4}
Let $q,q'\in \mH_0$. Then  $V$ is finite, nonnegative and satisfies
\[
\lb{deVS}
V={32\/\pi}\int_0^\iy zu(z)v^3(z)\,dz.
\]
\end{lemma}
\no{\bf Proof.} Since $q'\in \mH_0$, then  $H_2$ is finite and $I\in \ell_{3\/2}^1$. So $P_3$  is finite,  as well as
$V$. To prove \er{deVS} we start with a finite-gap aproximation for the momentum spectrum $\s_M=\R\sm \cup g_n$  of the potential $q_0+q$,
obtained by closing the gaps $g_m$ with large $m$. Namely, we fix $r>0$ and consider a new momentum spectrum $\s_M^r=\s_M\cup(-\iy,-r)\cup (r,\iy)$, where the new gaps $g_n^r$ are given by
$$
g_n^r=\ca g_n & if \ \  g_n\ss (-r,r) \\
        \es & if \ \ g_n\not\ss (-r,r)  \ac.\qqq 
$$
The variables corresponding to $\s_M^r$ will be indicated the upper index $r$.
Due to the general construction, presented in Section 2.2, for the finite-gap  momentum spectrum $\s_M^r$ there exists a unique conformal mapping
$$
k^r: \C\sm g^r\to \C\sm \G_n^r,\qqq  \G_n^r=(u_n^r+ih_n^r, u_n^r-ih_n^r), \qq h_n^r\ge 0,
$$
satisfying the asymptotics $k^r(z)=z-{O(1) z^{-1}}$ as $|z|\to \iy$. By \er{Kq},  each function  $k^r(z)-z$
is analytic
at $\iy$.  The  sequence of real numbers  $u_n^r, n\in \Z$ is odd, strongly increasing and $u_n^r\to\pm\iy$ as $n\to\pm\iy$.
In general, $k^r$ is not a quasimomentum for some periodic potential, since not necessarily $u_n=\pi n$ for all $n$.

For each $r$ we introduce  $Q_m^r, S_m^r, P_3^r$ and $V^r$ by relations \er{2.130},
\er{2.50} and \er{deVS} respectively, where $k=k^r, u=u^r$  and $v=v^r$.
Since $v^r(x)=0$ for large real $x$ and  $v^r(x+i0), u^r(x)\ge 0$ for real $x\ge 0$, then all these quantities
are finite and non-negative.
It is known (see \cite{L1}-\cite{L3}) that
$$
v^r(x) \nearrow v(x),\ \  |u^r(x)|\nearrow |u(x)| \qqq x\in \R, \qq \as \qq r\to\iy,
 $$
and that $k^r$ converges to $k$ uniformly on compact sets  from $\C\sm \s_M$.
From these convergence and Levy's theorem it follows that
\[
\lb{cQSPV}
Q_m^r\nearrow  Q_m,\qqq  S_m^r\nearrow  S_m,\qqq P_m^r\nearrow  P_m, \qqq  V^r \nearrow  V
\qq \as \qq r\to\iy,
\]
for $m=-1,0,1, 2,...$ (some limits may be infinite).

Assume that for each $r$ sufficiently large we have proved that
\[
\lb{x}
8(S_1^r-S_{-1}^rS_0^r)=P_3^r-V^r.
\]
Then sending $r\to \iy$ using \er{cQSPV} and evoking \er{H2},
we get that $H_2=P_3-($r.h.s. of \er{deVS}). Since $H_2=P_3-V$, we recover  \er{deVS}.

So it remains to show \er{x}.
Fix $r>1$ large enough and consider the integral $\int_{|z|=t}zk^4(z)\,dz$. The function $z(k^r(z))^4$ is analytic in
$\{|z|>r\}$. For any $m\ge 1$ we write its Tailor series at infinity,  omitting the index $r$ for brevity:
\[
\lb{km111} k(z)=z-{Q_{0}\/z}-{Q_{2}\/z^{3}}...-{Q_{2m}\/z^{2m+1}}+{O(1)\/z^{2m+2}} \qqq \as \qqq |z|\to \iy.
\]
Due to \er{km111} we get
$$
zk^4=z\rt(z^2-S_{-1}-{S_0\/z^2}-{S_1+o(1)\/z^4}\rt)^2=z^5\rt(1-{S_{-1}\/z^2}-{S_0\/z^4}-{S_1+o(1)\/z^6}\rt)^2
$$
$$
=z^5\rt(1-2\rt({S_{-1}\/z^2}+{S_0\/z^4}+{S_1+o(1)\/z^6}\rt)+\rt({S_{-1}\/z^2}+{S_0\/z^4}\rt)^2+..\rt)
=z^5+..-2{S_{1}-S_0S_{-1}\/z}+{O(1)\/z^2}.
$$
If $t>r$, then
\[
\lb{I1}
{1\/2\pi i}\int_{|z|=t}zk^4(z)\,dz=-2(S_{1}-S_0S_{-1}),
\]
since $v(z)=0$ for real  $z$ such that $|z|>r$. Thus
$$
{1\/2\pi i}\int_{|z|=t}zk^4(z)\,dz={1\/2\pi i}\int_{|z|=t}z(u^2-v^2+2iuv)^2\,dz={-1\/\pi i}\int_{\R}2z(u^2-v^2)(2iuv)dz
$$
\[
\lb{I2}
=-{8\/\pi }\int_{\cup_{n\ge 1} g_n} z(u^2-v^2)uv\,dz.
\]
By \er{H2}, \er{I1}, \er{2.50} and \er{I2}  we get that
$$
8(S_{1}-S_0S_{-1})={32\/\pi }\int_0^\iy z(u^2-v^2)uv\,dz={32\/\pi}\sum_{n\ge 1}\int_{g_n} z(u^2-v^2)uv\,dz=P_3-V,
$$
which yields \er{x}.
\BBox
\medskip

Let $0<a<{1\/4}$. We  note that since $|h\|_{2,a}^2=\sum_{n\ge 1}(2\pi n)^{2a-1}(2\pi n h_n^2)$, then
\[
\lb{hhn} \|h\|_{2,a}^2\le C_{2-4a} \|h\|_{4,1},\qqq where \qq
C_t^2=\sum _{n\ge 1}{1\/(2\pi n)^{t}}<\iy \qq {\rm if}\ t>1.
\]

\begin{theorem}
\lb{TVe}

A sequence   $h=h(I)$ belongs to  $\ell_1^4$  if and only if   $I\in \ell^2$.
 If $I\in \ell^2$, then the series
\[
\lb{Ven}
W=\sum_{n\ge 1}(4\pi n)V_n, 
\]
where $V_n=V_n(I) \ge 0$ is defined by \er{dVn}, converges for $I\in\ell^2_+$
 and defines in $\ell^2_+$ a finite non-negative function, equal to
$$
{32\/\pi}\int_0^\iy zu(z)v^3(z)\,dz.
$$
Moreover,

i)
the  function $W(I)$ satisfies the following estimates
\[
\lb{VhV}
{1\/5}\sum_{n\ge 1}(4\pi n)h_n^2I_n\le W\le 2\sum_{n\ge 1}(4\pi n)h_n^2I_n,
\]
\[
\lb{VA1}
W\le 4\|h\|_\iy^2P_1;
\]

ii) it  is continuous on $\ell_+^2$;

iii) on  the octant  $ \tilde \ell_+^\iy$ the function  $W(I)$  coincides with $V(I)$.

\end{theorem}
\no{\bf Proof.} Estimates \er{hIn} imply that $h\in \ell_1^4$ iff $I\in \ell^2$.

Due to the last inequality in   \eqref{VhVn},
$$
W
\le \sum_{n\ge 1}(4\pi n)2h_n^2I_n
\le 4\|I\|_2^{1/2}\|h\|_{4,1}^2,
$$
So $W$ is defined by a converging series  and satisfies the second estimate in \eqref{VhV}. Using the
lower bound for $V_n$ in  \eqref{VhVn} we recover the first estimate in \eqref{VhV}. Estimate \er{VA1} follows from \er{VhV}.

Since $u=\pi n$ on $g_n$ and $u$ vanishes outside $\cup g_n$, then $W(I)$ has the required integral representation.

 Let a sequence $I^s=(I_n^s)_1^\iy\xrightarrow{} I$ strongly in $\ell^2$ as $s\to \iy$. To prove ii) we  need to show that
\[
\lb{Vnhm} W(I^{s})\to W(I)\qqq \as \qqq s\to \iy.
\]
Using \er{hIn} and \er{3.0} we have
$$
\|h^s\|_{4,1}^4\le {3\pi^2\/8}C_0^{3\/2}\|I^s\|_2^2,
$$
 where $C_0\le \exp{\sqrt{2P_{-1}(I^s)}}$ and $P_{-1}(I)=\sum_{n\ge 1}{I_n}/{(2\pi n)}\le C_2\|I\|_2$.
Together with \er{hhn} this yields the estimates
\[
\lb{4.ex}
\sup _{s\ge 1}\|h^s\|_{4,1}<\iy,\qqq \sup _{s\ge
1}\|h^s\|_{2,a}<\iy \ \ {\rm if}\qq  a<{1\/4}.
\]
We claim that
\[
\lb{4.x}
h^{s}\xrightarrow{}  h \qqq {\rm  weakly \ in }\qq \ell_a^2 \qq \as \qq s\to\infty,
\]
for some $h\in \ell_a^2$. Indeed, assume that this is not the  case. Then by \er{4.ex} there are two different vectors $h', h''\in\ell_a^2$ and two subsequence $\{s_j'\}$ and $\{s_j''\}$ such that
\[
\lb{4.xc}
h^{s_{j}'}\xrightarrow{}  h', \ \ \ h^{s_{j}''}\xrightarrow{}  h''
 \qqq \qqq {\rm  weakly \ in }\qq \ell_a^2.
\]
Then
$$
h^{s_{j}'}\xrightarrow{}  h', \ \ \ h^{s_{j}''}\xrightarrow{}  h''
 \qqq \qqq {\rm  strongly\ in }\qq \ell_\n^2,
$$
for each $\n<a$. Using Theorem \ref{T1.1} and the identity
$$
k(z,h)=\int_0^z {\D'(t,h)\/\sqrt{1-\D^2(t,h)}}dt, \qqq z\in\cZ,
$$
which easily follows from \er{1.16}, we deduce that the corresponding conformal mappings  $k$ converge to limits:
\[
\lb{ccm}
k(z,h^{s_{j}'})\to k(z,h'), \qqq k(z,h^{s_{j}'})\to k(z,h')  \qqq\ as \ \ j\to\iy ,
\]
uniformly on bounded subsets in $\C$.
These convergences and \er{deac} imply that  for each $n$
the actions $I_n(h^{s_{j}'})$ and $I_n(h^{s_{j}''})$
converge to limits $I_n(h')$ and $I_n(h'')$, which must equal $I_n$.
That is, $h'$ and  $h''$ belong to the same iso-spectral class. Since  $h', h''\in \ell^2$, then by  Theorem \ref{T1.1} we have $h'=h''$. This proves \er{4.x}.

Due to \er{VnI}
$$
(4\pi n)V_n-(2\pi n)^3I_n={4\/\pi i}\int_{c_n}zk^4(z)\,dz=-{8\/\pi
i}\int_{\c_n}z^2(k,h)k^3dk.
$$
By this relation, \er{deac} and \er{ccm} we have
\[
\lb{4.xx}
V_n(I^{s_j})\to V_n(I)\qq  \as \qq j\to \iy \qq {\rm for \ each}  \qq n=1,2,3.....
\]
For any $N$ denote
$$
W^{(N)}(I^s)=\sum_{n\ge N}(4\pi n)V_n(I^s).
$$
Using \er{VhVn} we get
$$
W^{(N)}(I^s)\le \sum_{n\ge N}(8\pi n)(h_n^{s})^2 I_n^s=A_N+B_N,
$$
where
$$
A_N=\sum_{n\ge N} (8\pi n) (h_n^{s})^2 I_n,\qqq
B_N=\sum_{n\ge N}  (8\pi n)  (h_n^{s})^2 (I_n^s-I_n).
$$
Since
$$
A_N\le 4\|h^s\|_{4,1} \rt (\sum_{n>N}I_n^2\rt)^{1\/2},\qqq  B_N\le 4\|h^s\|_{4,1}\|I^s-I\|_2,
$$
then \er{4.ex} and \er{4.xx} yield the required convergence \eqref{Vnhm}.

The last assertion follows from \eqref{deVS} and the integral representation for $W(I)$.
 \BBox
 \medskip

{\bf Proof of Theorem \ref{T1}}.
  Theorem \ref{TVe} gives that the function $V: \wt\ell^\iy\to \R$ extends to a non-negative continuous function on the octant $\ell_+^{2}$.
Using estimates \er{VA1} and \er{Te1} we obtain $ V\le
4\|h\|_\iy^2P_1\le 8P_{-1}P_1, $ which yields \er{VA}.
If $I=0$, then \er{VA} implies $V(I)=0$. Finally, let $V=0$ for some $I$. Since the terms $V_n$ are non-negative, then each $V_n=0$ and \er{VhVn} implies that $I=0$.

It remains to prove \er{VA2}. Estimates \er{VhV} and \er{hIn} give
\[
\lb{21}
V\ge {2\/5}\sum_{n\ge 1}(2\pi n)h_n^2I_n  \ge {\pi\/10 C_I}\|I\|_2^2,
\]
which yields the first inequality in \er{VA2}.
Now we show the second. Using \er{VhV} and \er{Te5} we find that
$$
V\le \sum_{n\ge 1}(8\pi n h_n^2)I_n  \le \sum_{n\ge 1}
\rt(C_0^{1\/2}6\pi I_n^2+8{\|\r\|_2^2\/\pi}h_n^2I_n\rt)
\le 6\pi C_0^{1\/2}\|I\|_2^2+8{\|\r\|_2^2\/\pi}\|h\|_2\|h\|_\iy\|I\|_2.
$$
 Using \er{11}, \er{12} and   \er{Te1} we obtain
$$
\|\r\|_2^2\|h\|_2\|h\|_\iy\le \pi 4^4 \rt(1+Q_0^{1\/2}\rt) ^{1\/2}Q_0^2.
$$
Combining these  estimates we get that
$$
V\le 6\pi\sqrt{C_0}\|I\|_2^2+4^{11\/2}(1+Q_0^{1\/2})^{1\/2}Q_0^2\|I\|_2.
$$
Together with   \er{Te1} this yields the required estimate, since $C_0\le
C_-=\exp {\sqrt{2P_{-1}}}$.
\BBox
\medskip

Finally, as a by-product of some relations, derived above   in this work, we get two-sided algebraical
 bounds on the norm
$\|q'\|$ in terms of $P_3=\|I\|_{1,{3\/2}}$ (see \cite{K4} for two-sided algebraical estimates
of $\|q^{(m)}\|$ in terms of $P_{m+{1\/2}}$ for all $m\ge 0$).

\begin{proposition}
\lb{T3a}
The following estimates hold true:
\[
\lb{dqA}
\|q'\|^2\le 4(P_3+2P_1^2),
\]
\[
\lb{Adq}
P_3\le {\|q'\|^2\/2}
+{\| q'\|\/\sqrt 2}\|q\|^2
+2\pi \|q\|^3(1+\|q \|^{1\/3}).
\]

\end{proposition}

\no {\bf Proof}. Since $H_2=P_3-V$  and $\|q\|_\iy=\sup_{x\in [0,1]} |q(x)|\le {\| q'\|\/\sqrt 2}$, then
$$
{\|q'\|^2\/2}=H_2(q)-\int_0^1q^3(x)\,dx\le P_3+\|q\|_\iy\|q\|^2\le
P_3+{\| q'\|\/\sqrt 2}\|q\|^2\le P_3+{\|q\|^4\/2}+{\| q'\|^2\/4},
$$
which together with \er{pq}  yields \er{dqA}.
Using \er{VA1} and relations $\|q\|_\iy\le {\| q'\|\/\sqrt 2}$, $\|q\|^2=2P_1$ (see \er{pq}),  we obtain
$$
P_3={\|q'\|^2\/2}+\int_0^1q^3(x)\,dx-V\le {\|q'\|^2\/2}+{\| q'\|\/\sqrt 2}\|q\|^2+4\|h\|_\iy^2P_1=
{\|q'\|^2\/2}+{\| q'\|\/\sqrt 2}\|q\|^2+2\|h\|_\iy^2\|q\|^2.
$$
As
$\
\|h\|_\iy^2 \le  \pi \|q\|(1+\|q \|^{1\/3})
$  (see Theorem 2.3 in \cite{K3}), then we get \eqref{Adq}.
\BBox
\bigskip

\no  {\bf Acknowledgments.}\small Various parts of this paper were written at
CMLS Ecole Polytechnique,  France (April -- July, 2010).
E. K. is grateful CMLS  for their hospitality.
This work was supported by the
Ministry of education and science of the Russian Federation, state
contract 14.740.11.0581, and by l'Agence Nacionale de la Recherche, grant ANR-10-BLAN 0102.


\begin{thebibliography}
{999}\setlength{\itemsep}{-\parskip} \footnotesize


\bibitem[BoKu]{BoK} Bobenko, A.; Kuksin, S.  Finite-gap periodic solutions
of the KdV equation are nondegenerate. Phys. Lett. A  161  (1991),
no. 3, 274–-276.

\bibitem[BiKu]{BiK} Bikbaev, R.; Kuksin, S.
 On the parametrization of finite-gap solutions by frequency and wavenumber vectors and a theorem of I.
 Krichever. Lett. Math. Phys.  28  (1993),  no. 2, 115--122

\bibitem[FM]{FM} Flaschka H.; McLaughlin D. Canonically
conjugate variables for the Korteveg- de Vries equation and the Toda
lattice with periodic boundary conditions. Prog. of Theor. Phys.
55(1976),  438-456.




\bibitem[J] {J}  Jenkins, A. Univalent functions and conformal mapping.
Berlin, G\"ottingen, Heidelberg: Springer, 1958.


\bibitem[Ka] {Ka} Kappeler, T.
Fibration of the phase space for the Korteveg-de-Vries equation.
Ann. Inst. Fourier (Grenoble), 41, 1, 539-575 (1991).

\bibitem[KaP] {KP} Kappeler, T.; P\"oschel, J. KdV $\&$ KAM. Springer, 2003.

\bibitem[KaT]{KT} Kappeler, T.; Topalov, P.
Global wellposedness of KdV in $H^{-1}(\T,\R)$, Duke Math. J. 135
(2006), no. 2, 327--360.



\bibitem    [KK1] {KK1} Kargaev, P.; Korotyaev, E.  The inverse problem for the Hill operator, direct approach.
Invent. Math.  129, no. 3, 567-593(1997)

\bibitem[KK2]{KK2} Kargaev, P.; Korotyaev, E., Effective masses and
conformal mappings.  Commun. Math. Phys. 169(1995), 597-625.

\bibitem[K1]{K1} Korotyaev, E. Metric properties of conformal mappings  on the complex plane with parallel slits.
 Inter. Math. Reseach. Notices. 10(1996),  493-503.

 \bibitem[K2]{K2} Korotyaev E. Characterization of the spectrum of Schr\"odinger operators
with periodic distributions. Int. Math. Res. Not.  (2003) no. 37, 2019-2031.

\bibitem[K3]{K3} Korotyaev, E. Estimates for the Hill operator.I,
Journal Diff. Eq. 162(2000), 1-26.

\bibitem[K4]{K4}  Korotyaev, E.  Estimate for the Hill operator.II,
 J. Differential Equations 223 (2006), no. 2, 229-260.

\bibitem[K5]{K5} Korotyaev E. The estimates of  periodic potentials
in terms of effective masses. Commun. Math. Phys. 183(1997), 383-400.

\bibitem[K6]{K6}  Korotyaev, E. Estimate of periodic potentials in terms of gap lengths.
Commun. Math. Phys. 197(1998), no. 3, 521-526.

\bibitem[K7]{K7} Korotyaev, E. Periodic "weighted" operators. J.
Differential Equations 189 (2003), no. 2, 461--486.


\bibitem[Kr]{Kr} Krein M.: On the characteristic function $A(\l)$ of a linear canonical
system differential equation of the second order with periodic coefficients (Russian),
 Prikl. Mat. Meh. 21, 320-329 (1957).

\bibitem[Kri] {Kri}  Krichever, I. M. Perturbation theory in periodic problems for two-dimensional integrable
systems, Sov., Sci. Rev. C. Math. Phys. 9(1991), 1-101.



\bibitem[Ku] {Ku}  Kuksin, S. B. Analysis of Hamiltonian PDEs.
Oxford Lecture Series in Mathematics and its Applications, 19.
Oxford University Press, Oxford, 2000.

\bibitem[Ku1] {Ku1}  Kuksin, S. B. Damped-driven KdV and effective equations for long-time behaviour of its solutions. Geom. Funct. Anal.  20(2010),  no. 6, 1431-1463

\bibitem[KuP] {KuP}  Kuksin, S. B.; Perelman, G. Vey theorem in infinite dimensions and its application to
 {K}d{V}. DCDS-A 27 (2010), 1-24.




\bibitem[L1]{L1} Levin, B. Ya. Majorants in classes of subharmonic functions. (Russian) Teor. Funktsii Funktsional. Anal. i Prilozhen. No. 51 (1989), 3--17; translation in J. Soviet Math. 52 (1990), no. 6, 3441-3451.

\bibitem[L2]{L2} Levin, B. Ya.     The connection of a majorant with a conformal mapping. II.
(Russian) Teor. Funktsii Funktsional. Anal. i Prilozhen. No. 52 (1989), 3-21;
translation in J. Soviet Math. 52 (1990), no. 5, 3351--3364.

\bibitem[L3]{L3} Levin, B. Ya.
Classification of closed sets on $R$ and representation of a majorant. III. (Russian) Teor. Funktsii Funktsional. Anal. i Prilozhen. No. 52 (1989), 21--33; translation in J. Soviet Math. 52 (1990), no. 5, 3364--3372

\bibitem   [MO]{MO} Marchenko, V.; Ostrovski I. A characterization of the spectrum  of the Hill operator. Math. USSR
Sbornik  26(1975), 493--554.

\bibitem   [M]{M} Marchenko, V.
Sturm-Liouville operators and applications. Revised edition. AMS Chelsea Publishing, Providence, RI, 2011.



\bibitem[MT1]{MT1}  McKean H.; Trubowitz E.   Hill's operator and hyperelliptic function theory in the presence
of infinitely many branching points, Comm. Pure Appl. Math. 29, 1976, 143-226.


\bibitem[MT2]{MT}  McKean H.; Trubowitz E. Hill's surfaces and their theta
functions, Bull. Am. Math. Soc. 84, 1978, 1042-1085.

\bibitem[RS]{RS} Reed, M.;  Simon, B. Methods of modern mathematical
physics. IV. Analysis of operators, Academic Press, New York-London, 1978.

\bibitem[S]{S}
Stein, E. M. Functions of exponential type. Ann. of Math. (2)  65  (1957), 582–592.




\end{thebibliography}
\end{document}